\documentclass[11pt]{article}
\usepackage{amssymb, amsthm, amsmath}

\newtheorem{thm}{Theorem}[section]
\newtheorem{lem}[thm]{Lemma}

\newtheorem{cor}[thm]{Corollary}
\theoremstyle{definition}\newtheorem{df}[thm]{Definition}
\theoremstyle{definition}\newtheorem{rem}[thm]{Remark}
\theoremstyle{definition}\newtheorem{exm}[thm]{Example}

\renewcommand{\phi}{\varphi}

\newcommand{\N}{\mathbb{N}}
\newcommand{\Z}{\mathbb{Z}}
\newcommand{\Q}{\mathbb{Q}}
\newcommand{\R}{\mathbb{R}}
\newcommand{\T}{\mathbb{T}}

\newcommand{\Homeo}{\operatorname{Homeo}}
\newcommand{\Inf}{\operatorname{Inf}}
\newcommand{\Isom}{\operatorname{Isom}}
\newcommand{\id}{\operatorname{id}}

\newcommand{\Ker}{\operatorname{Ker}}

\newcommand{\xa}{(X,\alpha)}
\newcommand{\yb}{(Y,\beta)}

\title{Approximate conjugacy and full groups \\
of Cantor minimal systems
\thanks{2000 Mathematics Subject Classification: 37B05.}}
\author{MATUI Hiroki
\thanks{The author was supported by 
Grant-in-Aid for Young Scientists (B) of
Japan Society for the Promotion of Science.}}
\date{}

\begin{document}
\maketitle

\begin{abstract}
H. Lin and the author introduced the notion of approximate conjugacy 
of dynamical systems. 
In this paper, we will discuss the relationship between 
approximate conjugacy and full groups of Cantor minimal systems. 
An analogue of Glasner-Weiss's theorem will be shown. 
Approximate conjugacy of dynamical systems on the product 
of the Cantor set and the circle will also be studied. 
\end{abstract}

\section{Introduction}
In \cite{LM}, several versions of 
approximate conjugacy were introduced for minimal dynamical systems 
on compact metrizable spaces. 
In this paper, we will restrict our attention on dynamical systems 
on zero or one dimensional compact spaces and 
discuss approximate conjugacy. 

Let $X$ be the Cantor set. A homeomorphism $\alpha\in\Homeo(X)$ is 
said to be minimal when it has no nontrivial closed invariant sets. 
We call $\xa$ a Cantor minimal system. 
Giordano, Putnam and Skau introduced the notion of 
strong orbit equivalence for Cantor minimal systems in \cite{GPS1}, 
and showed that two systems are strong orbit equivalent if and only if 
the associated $K^0$-groups are isomorphic. 
We will show that this theorem can be regarded 
as an approximate version of Boyle-Tomiyama's theorem 
(\cite[Theorem 2.4]{GPS1} or \cite[Theorem 3.2]{BT}). 
Moreover, it will be also pointed out that 
two systems are strong orbit equivalent if and only if 
the closures of the topological full groups are isomorphic. 
This is an approximate analogue of \cite[Corollary 4.4]{GPS2}, 
in which it was proved that two systems are flip conjugate 
if and only if the topological full groups are isomorphic. 

In \cite{GW}, Glasner and Weiss proved that 
two Cantor minimal systems $\xa$ and $\yb$ are weakly orbit equivalent 
if and only if the associated $K^0$-groups are weakly isomorphic 
modulo infinitesimal subgroups. 
We will discuss relation of this result to approximate conjugacy. 
More precisely, it will be shown that there exists $\gamma\in[\alpha]$ 
which is conjugate to $\beta$ and whose associated orbit cocycle 
has at most one point of discontinuity 
if and only if there are $\sigma_n\in[[\alpha]]$ such that 
$\sigma_n\alpha\sigma_n^{-1}\rightarrow\beta$. 
Furthermore, this is also shown to be equivalent to the existence of 
a unital order surjection from $K^0\yb$ to $K^0\xa$. 

In the last section, dynamical systems on the product space of 
the Cantor set $X$ and the circle $\T$ are studied. 
H. Lin and the author investigated 
approximate conjugacy of minimal dynamical systems on the Cantor set 
or the circle in \cite{LM}. 
In the present paper, we will consider approximate conjugacy 
on their product space $X\times\T$. 
In the orientation preserving case, we will show that 
two systems are weakly approximately conjugate if and only if 
their periodic spectrum coincide. 
In the non-orientation preserving case, however, 
it is not enough to assume the same periodic spectrum 
to obtain weakly approximate conjugacy. 
A necessary and sufficient condition involving 
a $\Z_2$-extension of a Cantor minimal system will be given. 
We will continue to study such kind of dynamical systems and 
related crossed product $C^*$-algebras in \cite{M3}. 
Note that, however, our method is valid only for a skew product 
extension associated with a cocycle taking its values in $\Isom(\T)$. 
In general, every minimal homeomorphism on the product space of 
the Cantor set $X$ and the circle $\T$ is 
of the form $\alpha\times\phi$, where $\alpha\in\Homeo(X)$ is minimal 
and $\phi:X\rightarrow\Homeo(\T)$ is a (continuous) cocycle. 
We do not know when these kinds of dynamical systems are 
weakly approximately conjugate at present. 
\bigskip

\noindent
\textbf{Acknowledgment.}
I am grateful to Huaxin Lin from whom I have learned a lot 
through the collaboration \cite{LM}. 
I would like to thank Takeshi Katsura and Fumiaki Sugisaki 
for many helpful discussions.

\section{Preliminaries}

Let $X$ be a compact metrizable space. 
Equip $\Homeo(X)$ with the topology of pointwise convergence in norm 
on $C(X)$. Thus a sequence $\{\alpha_n\}_{n\in\N}$ in $\Homeo(X)$ 
converges to $\alpha$, if 
\[ \lim_{n\rightarrow\infty}
\sup_{x\in X}|f(\alpha_n^{-1}(x))-f(\alpha^{-1}(x))|=0 \]
for every complex valued continuous function $f\in C(X)$. 
This is equivalent to say that 
\[ \sup_{x\in X}d(\alpha_n(x),\alpha(x)) \]
tends to zero as $n\rightarrow\infty$, where $d(\cdot,\cdot)$ is 
a metric inducing the topology of $X$. 
When $X$ is the Cantor set, this is also equivalent to say that, 
for any clopen subset $U\subset X$, there exists $N\in\N$ such that 
$\alpha_n(U)=\alpha(U)$ for all $n\geq N$. 

\begin{df}[{\cite[Definition 3.1]{LM}}]
Let $\xa$ and $\yb$ be dynamical systems on compact metrizable 
spaces $X$ and $Y$. 
We say that $\xa$ and $\yb$ are weakly approximately conjugate, 
if there exist homeomorphisms $\sigma_n:X\rightarrow Y$ and 
$\tau_n:Y\rightarrow X$ such that 
$\sigma_n\alpha\sigma_n^{-1}$ converges to $\beta$ in $\Homeo(Y)$ 
and $\tau_n\beta\tau_n^{-1}$ converges to $\alpha$ in $\Homeo(X)$. 
\end{df}
\bigskip

Let us recall the definition of the $K^0$-group of 
a Cantor minimal system. 

\begin{df}
Let $\xa$ be a Cantor minimal system. We call 
\[ B_\alpha=\{f-f\alpha^{-1} \ : \ f\in C(X,\Z)\} \]
the coboundary subgroup and 
define the $K^0$-group of $\xa$ by 
\[ K^0\xa=C(X,\Z)/B_\alpha. \]
We write the equivalence class of $f\in C(X,\Z)$ in $K^0\xa$ by $[f]$, 
or $[f]_\alpha$ if we need to specify the minimal homeomorphism. 
\end{df}

The $K^0$-group is a unital ordered group 
equipped with the positive cone $K^0\xa^+$ and the order unit $[1_X]$. 

\begin{df}
Let $\xa$ be a Cantor minimal system. 
We define 
\[ \Inf(K^0\xa)=\{[f]\in K^0\xa:\mu(f)=0\mbox{ for every 
$\alpha$-invariant probability measure $\mu$} \} \]
and call it the infinitesimal subgroup. 
The quotient group $K^0\xa/\Inf(K^0\xa)$ is denoted 
by $K^0\xa/\Inf$ for short. 
\end{df}

In \cite{GPS1}, it was proved that 
$K^0\xa$ is a complete invariant for strong orbit equivalence and 
that $K^0\xa/\Inf$ is a complete invariant for orbit equivalence. 
\bigskip

We need to recall the idea of Kakutani-Rohlin partitions and 
Bratteli-Vershik models for Cantor minimal systems. 
The reader may refer to \cite{HPS} for the details. 
Let $\alpha\in\Homeo(X)$ be a minimal homeomorphism 
on the Cantor set $X$. 
A family of non-empty clopen subsets 
\[ {\cal P}=\{X(v,k):v\in V,1\leq k\leq h(v)\} \]
indexed by a finite set $V$ and natural numbers $k=1,2,\dots,h(v)$ 
is called a Kakutani-Rohlin partition, 
if the following conditions are satisfied: 
\begin{itemize} 
\item ${\cal P}$ is a partition of $X$. 
\item For all $v\in V$ and $k=1,2,\dots,h(v)-1$, we have 
$\alpha(X(v,k))=X(v,k+1)$. 
\end{itemize}
Let $R({\cal P})$ denote the clopen set $\bigcup_{v\in V}X(v,h(v))$ 
and call it the roof set of ${\cal P}$. 
For each $v\in V$, the family of clopen sets 
$X(v,1),X(v,2),\dots,X(v,h(v))$ is called a tower, 
and $h(v)$ is called the height of the tower. 
We may identify the label $v$ with the corresponding tower. 
One can divide a tower into a number of towers with the same height 
in order to obtain a finer partition. For example, if one needs to 
divide $E(v,1)$ into two clopen sets 
$O_1$ and $O_2=E(v,1)\setminus O_1$, then 
one can put $E(v_1,k)=\alpha^{k-1}(O_1)$ 
and $E(v_2,k)=\alpha^{k-1}(O_2)$ for every $k=1,2,\dots,h(v)$. 
We shall refer to this procedure as a division of a tower. 

Let $\{{\cal P}_n\}_{n\in\N}$ be a sequence of Kakutani-Rohlin 
partitions. We denote the set of towers in ${\cal P}_n$ by $V_n$ 
and clopen sets belonging to ${\cal P}_n$ by $E(n,v,k)$ ($v\in 
V_n,k=1,2,\dots,h(v)$). We say that $\{{\cal P}_n\}_{n\in\N}$ 
gives a Bratteli-Vershik model for $\alpha$, if the following are 
satisfied: 
\begin{itemize}
\item The roof sets $R({\cal P}_n)=\bigcup_{v\in V_n}E(n,v,h(v))$ 
form a decreasing sequence of clopen sets, 
which shrinks to a single point. 
\item ${\cal P}_{n+1}$ is finer than ${\cal P}_n$ for all $n\in\N$ 
as partitions, and $\bigcup_n{\cal P}_n$ generates the topology of $X$. 
\end{itemize}
Note that, by taking a subsequence of $\{{\cal P}_n\}_{n\in\N}$, 
we may further assume the following: 
\begin{itemize}
\item $R({\cal P}_{n+1})$ is contained in some $E(n,v,h(v))$ 
for all $n\in\N$. 
\end{itemize}
Therefore, when we put 
\[ \widetilde{\cal P}_n \ = \ \{E(n,v,k) \ : \ v\in V_n, \ 
1\leq k<h(v)\}\cup\{R({\cal P}_n)\}, \]
it is easily verified that $\widetilde{\cal P}_n$'s ($n=1,2,\dots$) 
also generate the topology. 
\bigskip

Let us recall the definition of (topological) full groups of 
Cantor minimal systems. 

\begin{df}[\cite{GPS2}]
Let $\xa$ be a Cantor minimal system. 
\begin{enumerate}
\item The full group $[\alpha]$ of $\xa$ is the subgroup of 
all homeomorphisms $\gamma\in\Homeo(X)$ that preserves 
every orbit of $\alpha$. 
To any $\gamma\in[\alpha]$ is associated a map $n:X\rightarrow \Z$, 
defined by $\gamma(x)=\alpha^{n(x)}(x)$ for $x\in X$. 
\item The topological full group $[[\alpha]]$ of $\xa$ is 
the subgroup of all homeomorphisms $\gamma\in[\alpha]$, 
whose associated map $n:X\rightarrow \Z$ is continuous. 
\end{enumerate}
\end{df}

In \cite{GPS2}, it was proved that $[\alpha]$ is 
a complete invariant for orbit equivalence and 
that $[[\alpha]]$ is a complete invariant for flip conjugacy. 

The following is a consequence of the Bratteli-Vershik model 
for $\xa$ (see \cite{HPS}). 

\begin{lem}\label{Hopfequiv}
Let $\xa$ be a Cantor minimal system. 
Let $U$ and $V$ be clopen subsets. 
\begin{enumerate}
\item $[1_U]=[1_V]$ in $K^0(X,\alpha )$ 
if and only if there exists 
$\gamma\in[[\alpha]]$ such that $\gamma(U)=V$. 
\item $[1_U]\leq[1_V]$ in $K^0(X,\alpha )$ 
if and only if there exists 
$\gamma\in[[\alpha]]$ such that $\gamma(U)\subset V$. 
\end{enumerate}
\end{lem}
\bigskip

We need the notion of the periodic spectrum of dynamical systems. 

\begin{df}
Let $X$ be a compact metrizable space and 
let $\alpha$ be a homeomorphism on $X$. 
By the periodic spectrum of $\xa$ or $\alpha$, we mean 
the set of natural numbers $p$ for which there are disjoint 
clopen sets $U,\alpha(U),\dots,\alpha^{p-1}(U)$ whose union 
is $X$. 
We denote the periodic spectrum of $\alpha$ by $PS(\alpha)$. 
\end{df}

When $\xa$ is a Cantor minimal system, 
it is well-known that $p\in PS(\alpha)$ if and only if 
$[1_X]$ is divisible by $p$ in $K^0\xa$. 
In fact, if $U,\alpha(U),\dots,\alpha^{p-1}(U)$ are 
disjoint clopen sets whose union is $X$, then $p[1_U]=[1_X]$. 
Conversely, if an integer valued continuous function $f$ 
satisfies $p[f]=[1_X]$, then there exists $g\in C(X,\Z)$ 
such that $pf-1_X=g-g\alpha^{-1}$, 
and the clopen subset $U=g^{-1}(p\Z)$ does the work.

\section{Approximate conjugacy for Cantor minimal systems}

We would like to discuss approximate conjugacy 
between Cantor minimal systems in this section. 
The following was shown in \cite{LM}. 
But we would like to present another proof which does not 
use any $C^*$-algebra theory. 
This theorem and its proof will be extended to 
dynamical systems on the product of the Cantor set and 
the circle in the next section. 

\begin{thm}\label{wac}
Let $\xa$ and $\yb$ be Cantor minimal systems. 
Then the following are equivalent. 
\begin{enumerate}
\item There exists a sequence of homeomorphisms 
$\sigma_n:X\rightarrow Y$ such that $\sigma_n\alpha\sigma_n^{-1}$ 
converges to $\beta$ in $\Homeo(Y)$. 
\item The periodic spectrum $PS(\beta)$ of $\beta$ is contained 
in the periodic spectrum $PS(\alpha)$ of $\alpha$. 
\end{enumerate}
\end{thm}
\begin{proof}
(1)$\Rightarrow$(2). 
Suppose $\sigma_n\alpha\sigma_n^{-1}\rightarrow\beta$. 
When $p$ belongs to $PS(\beta)$, 
there exists a clopen set $U\subset Y$ such that 
$U\cap \beta^k(U)$ is empty for all $k=1,2,\dots,p-1$ and 
$U\cup\beta(U)\cup\dots\cup\beta^{p-1}(U)=Y$. 
By the assumption, we can find a natural number $n$ such that 
\[ \sigma_n\alpha\sigma_n^{-1}(\beta^k(U))=\beta^{k+1}(U) \]
for all $k\in\N$. Put $V=\sigma_n^{-1}(U)$. 
Then $V,\alpha(V),\dots,\alpha^{p-1}(V)$ are mutually disjoint and 
their union is $X$. 
Hence $p$ is in the periodic spectrum $PS(\alpha)$. 

(2)$\Rightarrow$(1). 
Let ${\cal F}$ be a clopen partition of $Y$. 
It suffices to show that there exists a homeomorphism 
$\sigma:X\rightarrow Y$ such that $\sigma\alpha\sigma^{-1}(U)=\beta(U)$ 
for all $U\in{\cal F}$. 
We can find a Kakutani-Rohlin partition 
\[ {\cal Q}=\{Y(w,l):w\in W,l=1,2,\dots,h(w)\} \]
such that $\widetilde{\cal Q}$ is finer than ${\cal F}$, 
where 
\[ \widetilde{\cal Q}=\{Y(w,l):w\in W,l=1,2,\dots,h(w)-1\}
\cup \{R({\cal Q})\}. \]
We will show $\sigma\alpha\sigma^{-1}(U)=\beta(U)$ 
for all $U\in\widetilde{\cal Q}$. 
Let $p$ be the greatest common divisor of $h(w)$'s. 
Then $p$ is clearly in the periodic spectrum of $\beta$, 
and so of $\alpha$. 
Furthermore, there is $N\in\N$ such that 
\[ \{pn:n\geq N\}\subset
\left\{\sum_{w\in W}a_wh(w):a_w\in\N\right\}. \]

By choosing a sufficiently small roof set, we can find a Kakutani-Rohlin 
partition 
\[ {\cal P}=\{X(v,k):v\in V,k=1,2,\dots,h(v)\} \]
for $\xa$ such that $h(v)$ is divisible by $p$ and not less than $pN$ 
for all $v\in V$. By the choice of $N$, we have 
\[ h(v)=\sum_{w\in W}a_{v,w}h(w) \]
for some natural numbers $a_{v,w}$ ($v\in V,w\in W$). 
Put $b_w=\sum_{v\in V}a_{v,w}$ and 
divide each tower corresponding to $w\in W$ into 
$b_w$ towers. Let us denote the resulting Kakutani-Rohlin partition 
by ${\cal Q}'$. Then 
\[ \#{\cal P}=\sum_{v\in V}h(v)=\sum_{v\in V}\sum_{w\in W}a_{v,w}h(w)
=\sum_{w\in W}b_wh(w)=\#{\cal Q}'. \]
Therefore there exists a bijective map $\pi$ from ${\cal Q}'$ to 
${\cal P}$ such that all consecutive two clopen sets in any towers of 
${\cal Q}'$ go to consecutive clopen sets in a tower of ${\cal P}$. 
Since all non-empty clopen subsets of the Cantor set are homeomorphic, 
there is a homeomorphism $\sigma:X\rightarrow Y$ such that 
$\sigma(\pi(U))=U$ for all $U\in {\cal Q}'$. 
It is not hard to see $\sigma\alpha\sigma^{-1}(U)=\beta(U)$ 
for all $U\in\widetilde{Q}$. 
\end{proof}

\begin{cor}
Two Cantor minimal systems $\xa$ and $\yb$ are weakly approximately 
conjugate if and only if $PS(\alpha)=PS(\beta)$. 
\end{cor}

The theorem above says that the convergence of 
$\sigma_n\alpha\sigma_n^{-1}$ to $\beta$ does not bring 
any $K$-theoretical information except for periodic spectrum. 
To involve $K^0$-groups, we have to put an assumption 
on the conjugating map $\sigma_n$. 
\bigskip

The following lemma was proved in \cite{LM}. 

\begin{lem}\label{key}
Let $X$ be the Cantor set and $\alpha$, $\beta$ be 
minimal homeomorphisms. 
Let $\cal P$ be a clopen partition of $X$. 
If $[1_U]=[1_{\beta(U)}]$ in $K^0(X,\alpha)$ for all $U\in \cal P$, 
then one can find a homeomorphism $\sigma\in[[\alpha]]$ 
such that $\sigma\alpha\sigma^{-1}(U)=\beta(U)$ 
for all $U\in \cal P$. 
\end{lem}

As an immediate consequence of this lemma, we have the following. 
We would like to give a proof using only dynamical terminology, 
while a similar result can be found in \cite{LM}. 

\begin{thm}
For Cantor minimal systems $\xa$ and $\yb$, 
the following are equivalent. 
\begin{enumerate}
\item $K^0\xa$ is unital order isomorphic to $K^0\yb$. 
\item $\xa$ and $\yb$ are strong orbit equivalent. 
\item There is a homeomorphism $F:X\rightarrow Y$ such that 
$\overline{[[\alpha]]}=F^{-1}\overline{[[\beta]]}F$. 
\item There exist homeomorphisms $\sigma_n\in[[\alpha]]$, 
$\tau_n\in[[\beta]]$ and $F:X\rightarrow Y$ 
such that $F\sigma_n\alpha\sigma_n^{-1}F^{-1}\rightarrow\beta$ 
and $F^{-1}\tau_n\beta\tau_n^{-1}F\rightarrow\alpha$ 
as $n\rightarrow\infty$. 
\end{enumerate}
\end{thm}
\begin{proof}
(1)$\Leftrightarrow$(2) follows \cite[Theorem 2.1]{GPS1}. 

(2)$\Rightarrow$(3). 
Suppose that the strong orbit equivalence between $\xa$ and $\yb$ 
is implemented by a homeomorphism $F:X\rightarrow Y$. 
Then we have $\{g\circ F:g\in B_\beta\}=B_\alpha$ 
(see (i)$\Rightarrow$(ii) of \cite[Theorem 2.1]{GPS1} or 
(2)$\Rightarrow$(3) of Theorem \ref{analogue1}), 
which implies that, for every $\tau\in[[\beta]]$ and $U\subset X$, 
\[ 1_U-1_U\circ F^{-1}\tau F=
(1_{F(U)}-1_{F(U)}\circ\tau)\circ F \]
is in the coboundary subgroup $B_\alpha$. 
As mentioned in \cite[Proposition 2.11]{GPS2}, we know 
\[ \overline{[[\alpha]]}=\{\gamma\in\Homeo(X):
[1_U]_\alpha=[1_U\circ\gamma]_\alpha
\mbox{ for all clopen sets }U\subset X\}, \]
and so $F^{-1}\tau F$ belongs to $\overline{[[\alpha]]}$. 
The other inclusion follows in a similar fashion. 

(3)$\Rightarrow$(4). 
Let $\{{\cal P}_n\}$ be a sequence of clopen partitions of $Y$ 
which generates the topology of $Y$. 
By applying Lemma \ref{key} to the minimal homeomorphism 
$F^{-1}\beta F$ and $F^{-1}({\cal P}_n)$, 
we obtain $\sigma_n\in[[\alpha]]$ such that 
\[ \sigma_n\alpha\sigma_n^{-1}(F^{-1}(U))=F^{-1}\beta F(F^{-1}(U)) \]
for every $U\in{\cal P}_n$. 
Thus $F\sigma_n\alpha\sigma_n^{-1}F^{-1}(U)=\beta(U)$ 
for every $U\in{\cal P}_n$. 
The homeomorphisms $\tau_n$ can be constructed in the same way. 

(4)$\Rightarrow$(1). 
It is easy to see that $\{g\circ F:g\in B_\beta\}=B_\alpha$. 
Since $K^0\xa=C(X,\Z)/B_\alpha$ and $K^0\yb=C(Y,\Z)/B_\beta$, 
the assertion is clear. 
\end{proof}
The theorem above can be vied as an approximate version of 
\cite[Theorem 2.4]{GPS1}, in which it was shown that 
$\alpha$ and $\beta$ are flip conjugate if and only if 
$[[\alpha]]=F^{-1}[[\beta]]F$. 
Furthermore, in \cite{GPS2}, it was proved that 
if $[[\alpha]]$ is isomorphic to $[[\beta]]$ as an abstract group, 
then there exists a homeomorphism $F:X\rightarrow Y$ such that 
$[[\alpha]]=F^{-1}[[\beta]]F$. 
Indeed every isomorphism between $[[\alpha]]$ and $[[\beta]]$ is 
implemented by a homeomorphism. 
By using the same argument as in \cite{GPS2}, 
one can prove a similar result for the closure of 
the topological full group: 
every isomorphism between $\overline{[[\alpha]]}$ and 
$\overline{[[\beta]]}$ is implemented by a homeomorphism. 
Hence the conditions in the theorem above are also equivalent to 
$\overline{[[\alpha]]}$ being isomorphic to $\overline{[[\beta]]}$ 
as an abstract group. 
\bigskip

The following theorem is an analogue of \cite[Theorem 2.3 (a)]{GW}. 

\begin{thm}\label{analogue1}
When $\xa$ and $\yb$ are Cantor minimal systems, 
the following are equivalent. 
\begin{enumerate}
\item There is a unital order homomorphism $\rho$ 
from $K^0\yb$ to $K^0\xa$. 
\item There exist a continuous map $F:X\rightarrow Y$ and 
a minimal homeomorphism $\gamma\in[\alpha]$ such that 
the integer valued cocycle associated with $\gamma$ 
has at most one point of discontinuity and $\beta F=F\gamma$. 
\item There exists a continuous map $F:X\rightarrow Y$ such that 
$\{g\circ F:g\in B_\beta\}$ is contained in $B_\alpha$. 
\item There exists a sequence of homeomorphisms 
$\sigma_n:X\rightarrow Y$ such that 
$\sigma_n\alpha\sigma_n^{-1}\rightarrow\beta$ 
and $\lim_{n\rightarrow\infty}[f\circ\sigma_n]_\alpha$ 
exists for all $f\in C(Y,\Z)$. 
\end{enumerate}
\end{thm}
\begin{proof}
(1)$\Rightarrow$(2). 
Let $\rho$ be a unital order homomorphism from $K^0\yb$ to $K^0\xa$. 
The proof goes in a similar fashion to \cite[Proposition 2.9]{GW}. 
Take $a_0\in X$ and $b_0\in Y$ arbitrarily 
and put $a_1=\alpha(a_0),b_1=\beta(b_0)$. 
For each $n\in\N$, we would like to construct 
a Kakutani-Rohlin partition 
${\cal P}_n=\{Y(n,v,k):v\in V_n,1\leq k\leq h(v)\}$ for $\yb$ 
and homeomorphisms $F_n:X\rightarrow Y$, $\gamma_n\in[[\alpha]]$, 
so that the following conditions are satisfied. 
\begin{enumerate}
\item[(a)] $\{{\cal P}_n\}_n$ gives a Bratteli-Vershik model 
for $\yb$, 
and the roof sets $R({\cal P}_n)$ shrink to $\{b_0\}$. 
\item[(b)] $b_0$ and $b_1$ belong to distinct towers 
$v_n^0\in V_n$ and $v_n^1\in V_n$ in every ${\cal P}_n$. 
\item[(c)] For all $v\in V_n$ and $k=1,2,\dots,h(v)$, 
$\rho([1_{Y(n,v,k)}]_\beta)=[1_{Y(n,v,k)}\circ F_n]_\alpha$. 
\item[(d)] For all $v\in V_n$ and $k=1,2,\dots,h(v)-1$, 
$\gamma_nF_n^{-1}(Y(n,v,k))=F_n^{-1}(Y(n,v,k+1))$. 
\item[(e)] $F_n(a_0)\in Y(n,v_n^0,h(v_n^0))$ 
and $F_n(a_1)\in Y(n,v_n^1,1)$. 
\item[(f)] If $Y(n,v,k)$ contains $Y(n+1,w,j)$, then 
$F_n^{-1}(Y(n,v,k))$ also contains $F_{n+1}^{-1}(Y(n+1,w,j))$. 
\item[(g)] $\gamma_{n+1}(x)=\gamma_n(x)$ 
for all $x\in X\setminus F_n^{-1}(R({\cal P}_n))$. 
\item[(h)] The clopen sets $F_n^{-1}(R({\cal P}_n))$ shrink to $\{a_0\}$ 
and the clopen sets $\gamma_nF_n^{-1}(R({\cal P}_n))$ shrink to $\{a_1\}$. 
\item[(i)] $\{\gamma_n^{i-1}(a_1):1\leq i\leq h(v_n^1),n\in\N\}$ 
is dense in $X$. 
\end{enumerate}
If this is done, we can finish the proof as follows. 
Define $\gamma\in[\alpha]$ by $\gamma(x)=\gamma_n(x)$ 
for $x\in X\setminus F_n^{-1}(R({\cal P}_n))$ and $\gamma(a_0)=a_1$. 
Then $\gamma$ is a well-defined homeomorphism by (g) and (h). 
Moreover, the associated integer valued map is clearly continuous 
on $X\setminus \{a_0\}$. By (i), the orbit of $a_1$ by $\gamma$ 
is dense in $X$. 
By (a) and (f), $F_n$ converges to a continuous surjection $F$ 
satisfying $F^{-1}(Y(n,v,k))=F_n^{-1}(Y(n,v,k))$. 
By (h), $F^{-1}(b_0)=\{a_0\}$ and $F^{-1}(b_1)=\{a_1\}$. 
It follows that $\gamma$ is a minimal homeomorphism. 
By (d), we have $\beta F=F\gamma$. 

Let us construct ${\cal P}_n$, $F_n$ and $\gamma_n$ inductively. 
Suppose that ${\cal P}_n$, $F_n$ and $\gamma_n$ has been constructed. 
For each $v\in V_n$, $\{F_n^{-1}(Y(n,v,k)):1\leq j\leq h(v)\}$ is 
a tower with respect to $\gamma_n$, that is, we have 
$\gamma_nF_n^{-1}(Y(n,v,k))=F_n^{-1}(Y(n,v,k+1))$ for $k\neq h(v)$ 
by (d). But each level $F_n^{-1}(Y(n,v,k))$ may not so small. 
In order to achieve (h) and (i), 
we divide each tower by using $\gamma_n$ 
so that every level $F_n^{-1}(Y(n,v,k))$ becomes sufficiently small. 
Let ${\cal Q}$ be the obtained clopen partition of $X$ and 
let $c_v$ be the number of towers 
which the original tower corresponding to $v\in V_n$ is divided into. 
Put 
\[ \iota=\min\{\mu(O):O\in {\cal Q},
\mu\mbox{ is an $\alpha$-invariant probability measure}\}. \]
Take a Kakutani-Rohlin partition 
${\cal P}_{n+1}=\{Y(n+1,w,j):w\in V_{n+1},1\leq j\leq h(w)\}$ 
for $\yb$ so that the following are satisfied. 
\begin{itemize}
\item The roof set $R({\cal P}_{n+1})$ is 
a sufficiently small clopen neighborhood 
of $b_0$ and contained in $Y(n,v_n^0,h(v_n^0))\cap
\beta^{-1}(Y(n,v_n^1,1))$. 
\item $\nu(R({\cal P}_{n+1}))<\iota$ 
for every $\beta$-invariant probability measure $\nu$.
\item The tower $v_{n+1}^1\in V_{n+1}$ goes through 
every tower $v\in V_n$ at least $c_v$ times. 
\item ${\cal P}_{n+1}$ is sufficiently finer than 
${\cal P}_n$ as a partition and $v_{n+1}^0\neq v_{n+1}^1$. 
\end{itemize}
The first three conditions can be achieved 
by taking a sufficiently small roof set. 
The last condition is done by dividing towers. 
Let $O_0\in{\cal Q}$ (resp. $O_1\in{\cal Q}$) be the clopen set 
that contains $a_0$ (resp. $a_1$). 
Since $\rho([1_{R({\cal P}_{n+1})}]_\beta)<[O_0]_\alpha,[O_1]_\alpha$, 
we can find a clopen neighborhood $O_0'$ (resp. $O_1'$) of 
$a_0$ (resp. $a_1$) which is contained in $O_0$ (resp. $O_1$) 
and whose $K^0$-class is equal to $\rho([1_{R({\cal P}_{n+1})}]_\beta)$. 
We will define $F_{n+1}$ so that $F_{n+1}(O_0')=R({\cal P}_{n+1})$ and 
$F_{n+1}(O_1')=\beta(R({\cal P}_{n+1}))$. 
By using Lemma \ref{Hopfequiv}, 
take $\sigma\in[[\alpha]]$ with $\sigma(O_0')=O_1'$, and define 
$\gamma_{n+1}$ on $X\setminus (R({\cal P}_n)\setminus O_0')$ by 
$\gamma_{n+1}(x)=\gamma_n(x)$ for $x\in X\setminus R({\cal P}_n)$ 
and $\gamma_{n+1}(x)=\sigma(x)$ for $x\in O_0'$. 
We need not be careful about the choice of $\sigma$, 
because in the next step we will replace $\gamma_{n+1}|O_0'$ 
by another one. 
Let us consider the tower $v_{n+1}^1$. 
By repeating use of Lemma \ref{Hopfequiv}, 
we can take a copy of the tower via $\rho$, and 
define $F_{n+1}$ and $\gamma_{n+1}$ on it so that (c), (d), (e) 
and (f) are achieved. 
Moreover, we can do it so that the tower goes through 
every tower in ${\cal Q}$, that is, each clopen set of ${\cal Q}$ 
intersects with $\{\gamma_{n+1}^{i-1}(a_1):1\leq i\leq h(v_{n+1}^1)\}$. 
Hence we can ensure the condition (i). 
Thereby the induction step is completed. 

(2)$\Rightarrow$(3). 
Suppose that the integer valued function $n:X\rightarrow\Z$ 
associated with $\gamma$ is continuous on $X\setminus\{a_0\}$. 
It suffices to show that, for every clopen set $U\subset Y$, 
$(1_U-1_{\beta(U)})\circ F$ belongs to $B_\alpha$. 
Put $V=F^{-1}(U)$. Then we have 
\[ F^{-1}\beta(U)=\gamma F^{-1}(U)=
\bigcup_{k\in\Z}\alpha^k(V\cap n^{-1}(k)). \]
If $V$ does not contain $a_0$, 
in the right-hand side of the equality above, the union is actually finite. 
Hence we can see that 
\[ (1_U-1_{\beta(U)})\circ F=
1_V-\sum_{k\in\Z}1_{V\cap n^{-1}(k)}\circ\alpha^{-k} \]
is in the coboundary subgroup $B_\alpha$. 
In the case that $V$ contains $a_0$, by means of 
\[ (1_U-1_{\beta(U)})\circ F=-(1_{U^c}-1_{\beta(U^c)})\circ F, \]
we can apply the same argument. 

(3)$\Rightarrow$(4). 
Let ${\cal P}_n$ be a sequence of clopen partitions of $Y$ such that 
${\cal P}_{n+1}$ is finer than ${\cal P}_n$ for all $n\in\N$ and 
$\bigcup{\cal P}_n$ generates the topology of $Y$. 
Let ${\cal Q}_n$ be the clopen partition generated by ${\cal P}_n$ 
and $\beta({\cal P}_n)$. 
Although $F$ is not a homeomorphism, for each $n\in\N$, 
we can find a homeomorphism 
$F_n:X\rightarrow Y$ such that $F_n^{-1}(U)=F^{-1}(U)$ 
for every $U\in{\cal Q}_n$. 
By the assumption, 
\[ 1_{F_n^{-1}(U)}-1_{F_n^{-1}\beta(U)}=
1_{F^{-1}(U)}-1_{F^{-1}\beta(U)}=(1_U-1_{\beta(U)})\circ F \]
is zero in $K^0\xa$ for all $U\in{\cal Q}_n$. 
Now Lemma \ref{key} applies to 
the minimal homeomorphism $F_n^{-1}\beta F_n$ and 
the clopen partition $F_n^{-1}({\cal Q}_n)$ of $X$ and 
yields $\tau_n\in[[\alpha]]$ such that 
\[ \tau_n\alpha\tau_n^{-1}(F_n^{-1}(U))=
F_n^{-1}\beta F_n(F_n^{-1}(U)) \]
for all $U\in{\cal Q}_n$. 
Put $\sigma_n=F_n\tau_n$. 
Then we get 
\[ \sigma_n\alpha\sigma_n^{-1}(U)=\beta(U) \]
for all $U\in{\cal Q}_n$, and moreover 
\[ [1_U\circ\sigma_n]_\alpha=[1_U\circ F_n\circ\tau_n]_\alpha=
[1_U\circ F_n]_\alpha=[1_U\circ F]_\alpha, \]
which completes the proof. 

(4)$\Rightarrow$(1). 
The unital order homomorphism $\rho$ is given by 
$\rho([g]_\beta)=\lim_{n\rightarrow\infty}[g\circ\sigma_n]_\alpha$ 
for $[g]\in K^0\yb$. 
By $\sigma_n\alpha\sigma_n^{-1}\rightarrow\beta$, 
one checks that this is really a homomorphism 
from $K^0\yb$. 
\end{proof}

\begin{rem}\label{elliott}
It is known that (1) directly implies (3) in the theorem above. 
See \cite[Theorem 2.6]{LM} for example. 
This fact can be understood as a corollary of Elliott's classification 
theorem of AF algebras in the following way. 
Let $A$ and $B$ be the unital AF algebras 
such that $K_0(A)$ and $K_0(B)$ are unital order isomorphic to 
$K^0\xa$ and $K^0\yb$, respectively. 
One can regard $C(X)$ and $C(Y)$ 
as the diagonal subalgebras of $A$ and $B$. 
Then, by the classification theorem, there is an homomorphism 
$\phi:B\rightarrow A$ such that $\phi_*=\rho$ and $\phi(C(Y))=C(X)$. 
Hence there exists a continuous map $F:X\rightarrow Y$ such that 
$\phi(g)=g\circ F$ for all $g\in C(Y)$. 
It is not hard to see that $F$ meets the requirement. 
\end{rem}
\bigskip

Next, we would like to consider 
when one can choose $F$ to be a homeomorphism 
in the conditions (2) and (3) in Theorem \ref{analogue1}. 

\begin{df}
Let $G$ and $H$ be dimension groups. 
We say an order homomorphism $\rho:G\rightarrow H$ is 
an order surjection, 
if for every $g\in G$ and $h\in H$ with $0\leq h\leq \rho(g)$ 
there is $g'\in G$ such that $\rho(g')=h$ and $0\leq g'\leq g$. 
\end{df}
Evidently the order surjectivity implies $\rho(G_+)=H_+$. 
The other implication, however, does not hold in general. 

\begin{exm}
Let $G$ be the subgroup of $C([0,\pi],\R)$ generated 
by $\Q[x]$ (the set of rational polynomials) and $2(\pi-x)$. 
Thus every $g(x)\in G$ has the form $f(x)+2n(\pi-x)$ 
with $f(x)\in\Q[x]$ and $n\in\Z$. 
Put 
\[ G^+=\{g(x)\in G:g(x)>0\mbox{ for all }x\in[0,\pi]\}\cup
\{0\}. \]
Then one checks that $(G,G^+,1)$ is a unital ordered group 
satisfying the Riesz interpolation property. 
Moreover $G$ is simple and has no infinitesimal elements. 

Put $H=\{g(\pi)\in\R:g(x)\in G\}$. 
We regard $H$ as a unital ordered group with the order 
inherited from $\R$. 
Then the point evaluation at $x=\pi$ gives a homomorphism 
$\rho$ from $G$ to $H$. 
It is obvious that $\rho$ is a unital order homomorphism. 
Clearly $\rho$ is surjective and its kernel is 
$\{2n(\pi-x):n\in\Z\}\cong\Z$, because $\pi$ is transcendental. 

Suppose that $h\in H$ is positive. Since $\rho$ is surjective, 
there is $g(x)\in G$ such that $\rho(g)=g(\pi)=h>0$. 
Although the function $g$ may not be positive, 
for a sufficiently large $n\in\N$, $g(x)+2n(\pi-x)$ is strictly 
positive on $[0,\pi]$. Hence we have $\rho(G^+)=H^+$. 
Nevertheless the order homomorphism $\rho$ is not 
an order surjection. 
To explain it, let $g(x)=x+1\in G$ and $h=4-\pi\in H$. 
Then $0<h<\rho(g)=\pi+1$. 
But there are no integer $n$ such that 
$0<4-x+2n(\pi-x)<x+1$ for all $x\in[0,\pi]$. 
\end{exm}

\begin{thm}\label{analogue2}
Let $\xa$ and $\yb$ be Cantor minimal systems. 
Then the following are equivalent. 
\begin{enumerate}
\item There is a unital order surjection $\rho$ from $K^0\yb$ 
to $K^0\xa$. 
\item There is $\gamma\in [\alpha]$ such that the integer valued 
cocycle associated with $\gamma$ has at most one point of 
discontinuity and $\gamma$ is conjugate to $\beta$. 
\item There is a homeomorphism $F:X\rightarrow Y$ such that 
$\{g\circ F:g\in B_\beta\}$ is contained in $B_\alpha$. 
\item There are homeomorphisms $\sigma_n\in[[\alpha]]$ and 
$F:X\rightarrow Y$ such that 
$F\sigma_n\alpha\sigma_n^{-1}F^{-1}\rightarrow\beta$. 
\end{enumerate}
\end{thm}
\begin{proof}
(1)$\Rightarrow$(2). 
This is done by changing a part of the proof of (1)$\Rightarrow$(2) 
in Theorem \ref{analogue1}. 
We follow the notation used there. 
In order to make $F:X\rightarrow Y$ a homeomorphism, 
we have to require that 
the partition $F_n^{-1}({\cal P}_n)$ of $X$ is sufficiently 
finer in each inductive step. 
We will achieve it by dividing the towers and changing $F_n$. 
Let us focus our attention on a tower corresponding to $v\in V_n$. 
The clopen sets $F_n^{-1}(Y(n,v,k))$ ($k=1,2,\dots,h(v)$) 
are not so small at first, and so we must divide these clopen sets. 
For simplicity, suppose that $F_n^{-1}(Y(n,v,k))$ is divided into 
two clopen sets $X(n,v_1,k)$ and $X(n,v_2,k)$ so that 
\[ \gamma_n(X(n,v_1,k))=X(n,v_1,k+1), \ 
\gamma_n(X(n,v_2,k))=X(n,v_2,k+1) \]
for all $k=1,2,\dots,h(v)-1$. 
By the order surjectivity of $\rho$, we can find the clopen subsets 
$Y(n,v_1,k)$ and $Y(n,v_2,k)$ of $Y(n,v,k)$ such that 
\[ \beta(Y(n,v_1,k))=Y(n,v_1,k+1), \ 
\beta(Y(n,v_2,k))=Y(n,v_2,k+1) \]
for all $k=1,2,\dots,h(v)-1$ and 
\[ \rho([1_{Y(n,v_1,k)}]_\beta)=[1_{X(n,v_1,k)}]_\alpha, \ 
\rho([1_{Y(n,v_2,k)}]_\beta)=[1_{X(n,v_2,k)}]_\alpha \]
for all $k=1,2,\dots,h(v)$. 
Then we rearrange $F_n$ so that $F_n(X(n,v_i,k))=Y(n,v_i,k)$ 
for $i=1,2$ and $k=1,2,\dots,h(v)$. 
In this way, the tower in $X$ can be divided. 
Note that we need to take care of the condition (e), 
when dealing with the towers corresponding to $v_n^0$ and $v_n^1$. 

(2)$\Rightarrow$(3), (3)$\Rightarrow$(4) and (4)$\Rightarrow$(1) 
can be proved in the same way as in Theorem \ref{analogue1}.
\end{proof}

\begin{rem}
Of course, (3) of the theorem above is equivalent to 
$F^{-1}\overline{[[\beta]]}F\subset\overline{[[\alpha]]}$. 
\end{rem}

\begin{rem}
One can prove (1)$\Rightarrow$(3) directly in a similar fashion 
to \cite[Theorem 2.6]{LM}. 
See also Remark \ref{elliott}. 
\end{rem}

In \cite{GW}, so to say, a `modulo infinitesimal' version of 
the theorem above was discussed. 
More precisely, Glasner and Weiss showed that 
if there is a unital order homomorphism 
from $K^0\yb/\Inf$ to $K^0\xa/\Inf$, then there exists 
a minimal homeomorphism $\gamma\in[\alpha]$ such that 
the Cantor minimal system $(X,\gamma)$ admits $\yb$ as a factor. 
Besides, they proved that If $K^0\yb/\Inf$ and $K^0\xa/\Inf$ are 
unital order isomorphic, then $\gamma$ can be chosen 
to be conjugate to $\beta$. 

We can modify the argument in the proof of Theorem \ref{analogue1} 
and \ref{analogue2} 
so that it applies to the `modulo infinitesimal' case. 
As a consequence, we get the following. 
\begin{thm}
Let $\xa$ and $\yb$ be Cantor minimal systems. 
The following are equivalent. 
\begin{enumerate}
\item There is a unital order surjection from $K^0\yb/\Inf$ 
to $K^0\xa/\Inf$. 
\item There exists $\gamma\in[\alpha]$ such that $(X,\gamma)$ 
is conjugate to $\yb$. 
\end{enumerate}
\end{thm}

\section{Approximate conjugacy on $X\times\T$}

In this section, we will extend Theorem \ref{wac} to dynamical systems 
on the product space of the Cantor set $X$ and the circle. 
The crossed product $C^*$-algebra arising from 
this kind of dynamical system will be discussed in \cite{M3}. 
We identify the circle with $\T\cong\R/\Z$ and 
denote the distance from $t\in\T$ to zero by $|t|$. 
The finite cyclic group of order $m$ is denoted 
by $\Z_m\cong\Z/m\Z$ and 
may be identified with $\{0,1,\dots,m-1\}$. 

Define $o:\Homeo(\T)\rightarrow\Z_2$ by 
\[ o(\phi)=\begin{cases}
0 & \phi\mbox{ is orientation preserving}\\
1 & \phi\mbox{ is orientation reversing.} \end{cases} \]
Then the map $o(\cdot)$ is a homomorphism. 
Let $R_t$ denote the translation on $\T=\R/\Z$ by $t\in\T$. 
By $\Isom(\T)$ we mean the set of isometric homeomorphisms on $\T$. 
Thus, 
\[ \Isom(\T)=\{R_t:t\in\T\}\cup\{R_t\lambda:t\in\T\}, \]
where $\lambda\in\Homeo(\T)$ is defined by $\lambda(t)=-t$, 
and so $\Isom(\T)$ is isomorphic to the semidirect product of 
$\T$ by $\Z_2$. 
Equip $\Isom(\T)$ with the topology induced from $\Homeo(\T)$. 
Then $\Isom(\T)$ is a disjoint union of two copies of $\T$ 
as a topological space. 

\begin{df}
Let $\xa$ be a Cantor minimal system. 
Suppose that a continuous map $X\ni x\mapsto\phi_x\in\Isom(\T)$ 
is given. 
We denote the homeomorphism 
$(x,t)\mapsto(\alpha(x),\phi_x(t))$ on $X\times\T$ 
by $\alpha\times\phi$, 
and call $(X\times\T,\alpha\times\phi)$ 
the skew product extension of $\xa$ by $\phi$. 
\end{df}

When $\phi:X\rightarrow\Isom(\T)$ is a continuous map, 
the composition of $\phi$ and $o$ gives a continuous function 
from $X$ to $\Z_2$. 
We denote this $\Z_2$-valued function by $o(\phi)$. 
Under the identification of 
\[ C(X,\Z_2)/\{f-f\alpha^{-1}:f\in C(X,\Z_2)\} \]
with $K^0\xa/2K^0\xa$, 
an element of $K^0\xa/2K^0\xa$ is obtained from $o(\phi)$. 
We write it by $[o(\phi)]$ or $[o(\phi)]_\alpha$. 
If $o(\phi)(x)=0$ for all $x\in X$, there exists a continuous 
function $\xi:X\rightarrow\T$ such that $\phi_x=R_{\xi(x)}$ 
for every $x\in X$. 
In this case, let us denote the induced homeomorphism on $X\times\T$ 
by $\alpha\times R_\xi$ instead of $\alpha\times\phi$. 

\begin{df}
Let $\xa$ be a Cantor minimal system 
and let $\phi:X\rightarrow\Isom(\T)$ be a continuous map. 
We say that $\alpha\times\phi$ (or $\phi$) is orientation preserving 
when $[o(\phi)]$ is zero in $K^0\xa/2K^0\xa$. 
\end{df}

\begin{lem}\label{orientation}
In the above setting, 
suppose that $\alpha\times\phi$ is orientation preserving. 
Then there exists a continuous map $\xi:X\rightarrow\T$ such that 
$\alpha\times\phi$ is conjugate to $\alpha\times R_\xi$. 
\end{lem}
\begin{proof}
Since $[o(\phi)]$ is zero in $K^0\xa/2K^0\xa$, 
there exists a continuous function $f:X\rightarrow\Z_2$ 
such that $o(\phi)(x)=f(x)-f(\alpha(x))$ for all $x\in X$. 
Define a continuous map $\psi:X\rightarrow\Isom(\T)$ by 
\[ \psi_x=\begin{cases}
\id & f(x)=0 \\
\lambda & f(x)=1. \end{cases} \]
Then 
\[ o(\psi_{\alpha(x)}\phi_x\psi_x^{-1})=
f(\alpha(x))+o(\phi_x)-f(x)=0, \]
and so there exists $\xi\in C(X,\T)$ such that 
$\psi_{\alpha(x)}\phi_x=R_{\xi(x)}\psi_x$ for all $x\in X$. 
Thus, $\id\times\psi$ gives a conjugacy 
between $\alpha\times\phi$ and $\alpha\times R_\xi$. 
\end{proof}

Although the following theorem is actually contained 
in Theorem \ref{wacXT}, 
we would like to present it as a prototype. 
Note that every clopen subset of $X\times\T$ is of the form $U\times\T$ 
with a clopen set $U\subset X$. 
Hence the periodic spectrum $PS(\alpha\times\phi)$ agrees 
with $PS(\alpha)$. 

\begin{thm}\label{proto}
Let $\xa$ and $\yb$ be Cantor minimal systems and 
let $\phi:X\rightarrow\Isom(\T)$ and $\psi:Y\rightarrow\Isom(\T)$ 
be continuous maps. 
If both $\alpha\times\phi$ and $\beta\times\psi$ are 
orientation preserving, then the following are equivalent. 
\begin{enumerate}
\item There exist homeomorphisms $\sigma_n$ from $X\times\T$ to 
$Y\times\T$ such that $\sigma_n(\alpha\times\phi)\sigma_n^{-1}$ 
converges to $\beta\times\psi$ in $\Homeo(\T)$.
\item The periodic spectrum $PS(\beta)$ of $\beta$ is contained 
in the periodic spectrum $PS(\alpha)$ of $\alpha$. 
\end{enumerate}
\end{thm}
\begin{proof}
(1)$\Rightarrow$(2). 
The proof is the same as that of Theorem \ref{wac}, 
because of $PS(\alpha\times\phi)=PS(\alpha)$ and 
$PS(\beta\times\psi)=PS(\beta)$. 

(2)$\Rightarrow$(1). 
Thanks to Lemma \ref{orientation}, we may assume that there exist 
$\xi\in C(X,\T)$ and $\zeta\in C(Y,\T)$ such that 
$\phi_x=R_{\xi(x)}$ and $\psi_y=R_{\zeta(y)}$ for all $x\in X$ 
and $y\in Y$. 

Let ${\cal F}$ be a clopen partition of $Y$ and 
let $\varepsilon>0$. 
It suffices to find a homeomorphism $\sigma:X\rightarrow Y$ 
and a continuous function $\eta:X\rightarrow \T$ such that 
\[ \sigma\alpha\sigma^{-1}(U)=\beta(U) \]
for all $U\in{\cal F}$ and 
\[ |(\xi-\zeta\sigma)(x)-(\eta-\eta\alpha)(x)|<\varepsilon \]
for all $x\in X$. 
If such $\sigma$ and $\eta$ are found, 
then $\sigma\times R_\eta$ does the work. 
In fact, for $(x,t)\in X\times\T$, 
\begin{align*}
(\sigma\times R_\eta)(\alpha\times R_\xi)(x,t)
&=(\sigma\times R_\eta)(\alpha(x),t+\xi(x))\\
&=(\sigma(\alpha(x)),t+\xi(x)+\eta(\alpha(x)))
\end{align*}
is close to 
\begin{align*}
(\beta\times R_\zeta)(\sigma\times R_\eta)(x,t)
&=(\beta\times R_\zeta)(\sigma(x),t+\eta(x))\\
&=(\beta(\sigma(x)),t+\eta(x)+\zeta(\sigma(x))). 
\end{align*}

Let 
\[ {\cal Q}=\{Y(v,k):v\in V,k=1,2,\dots,h(v)\} \]
be a Kakutani-Rohlin partition for $\yb$ such that 
$h(v)>\varepsilon^{-1}$ for all $v\in V$ and 
\[ \widetilde{\cal Q}=
\{Y(v,k):v\in V,k=1,2,\dots,h(v)-1\}\cup\{R({\cal Q})\} \]
is finer than ${\cal F}$. 
By Theorem \ref{wac}, there exists a homeomorphism 
$\sigma:X\rightarrow Y$ such that 
\[ \sigma\alpha\sigma^{-1}(U)=\beta(U) \]
holds for all $U\in\tilde{\cal Q}$. 
Put $X(v,k)=\sigma^{-1}(Y(v,k))$. Then 
\[ {\cal P}=\{X(v,k):v\in V,k=1,2,\dots,h(v)\} \]
is evidently a Kakutani-Rohlin partition for $\xa$. 
Define a continuous function $\kappa$ 
from $\alpha(R({\cal P}))$ to $\T$ by 
\[ \kappa(x)=
\sum_{i=0}^{h(v)-1}(\zeta\sigma-\xi)(\alpha^i(x)) \]
for all $x\in X(v,1)$. 
Since $X$ is totally disconnected, 
there exists a continuous function $\tilde{\kappa}$ from 
$\alpha(R({\cal P}))$ to $\R$ such that 
\[ \tilde{\kappa}(x)+\Z=\kappa(x) \mbox{ and }
-1<\tilde{\kappa}(x)<1 \]
for all $x\in \alpha(R({\cal P}))$. 
We can define the continuous function $\eta:X\rightarrow\T$ 
by $\eta(x)=0$ for $x\in \alpha(R({\cal P}))$ and 
\[ \eta(\alpha^j(x))=
\sum_{i=0}^{j-1}(\zeta\sigma-\xi)(\alpha^i(x))
-\frac{j}{h(v)}\tilde{\kappa}(x)+\Z \]
for $x\in X(v,1)$ and $j=1,2,\dots,h(v)-1$. 
One can check that 
$|(\xi-\zeta\sigma)(x)-(\eta-\eta\alpha)(x)|$ is less than 
$\varepsilon$ for all $x\in X$, 
because $|h(v)^{-1}\tilde{\kappa}(x)|$ is less than $\varepsilon$. 
The proof is completed. 
\end{proof}
\bigskip

We would like to consider the general case. 
Let $\alpha\times\phi$ be as above. 
Define a homeomorphism on $X\times\Z_2$ by 
\[ \alpha\times o(\phi):(x,k)\mapsto(\alpha(x),k+o(\phi)(x)). \]
Then, if the $\Z_2$-valued continuous function $o(\phi)$ is not 
zero in $K^0\xa/2K^0\xa$, by \cite[Lemma 3.6]{M1}, 
$\alpha\times o(\phi)$ is a minimal homeomorphism 
on the Cantor set $X\times\Z_2$. 
The projection $\pi$ from $X\times\Z_2$ to the first coordinate $X$ 
gives a factor map. 
It is well known that $\pi$ induces a unital order embedding $\pi^*$ 
from $K^0\xa$ to $K^0(X\times\Z_2,\alpha\times o(\phi))$. 
In particular, $PS(\alpha)$ is contained in $PS(\alpha\times o(\phi))$. 

In general, if a continuous $\Z_m$-valued function $c:X\rightarrow\Z_m$ 
is given, then one can define $\alpha\times c$ on $X\times\Z_m$ 
in a similar fashion. 
This kind of dynamical system was studied in \cite{M1} and \cite{M2}. 

Although we need the following lemma in the case $m=2$, 
we would like to describe a general version. 
\begin{lem}\label{torsion}
Let $\xa$ be a Cantor minimal system and 
let $c:X\rightarrow\Z_m$ be a continuous function. 
Suppose that $\alpha\times c$ is a minimal homeomorphism 
on $X\times\Z_m$. Then we have 
\[ T(K^0(X\times\Z_m,\alpha\times c)/
\pi^*(K^0\xa))\cong\Z_m, \]
where $T(\cdot)$ means the torsion subgroup. 
Moreover, its generator is given by 
the $\Z$-valued continuous function 
\[ f_0(x,k)=\begin{cases}
1 & c(\alpha^{-1}(x))\neq0\mbox{ and }
k\in\{0,1,\dots,c(\alpha^{-1}(x))-1\} \\
0 & \mbox{otherwise}. \end{cases} \]
\end{lem}
\begin{proof}
Let us denote the homeomorphism $(x,k)\mapsto(x,k+1)$ by $\gamma$. 
It is clear that $\gamma$ commutes with $\alpha\times c$. 

Suppose that $f\in C(X\times \Z_m,\Z)$ gives a torsion element 
in $K^0(X\times\Z_m,\alpha\times c)/\pi^*(K^0\xa)$. 
There exist $n\in\N$ and $g\in C(X,\Z)$ such that 
\[ nf-g\circ\pi\in B_{\alpha\times c}, \]
which implies that 
\[ nf\circ\gamma^{-1}-g\circ\pi\circ\gamma^{-1}
=nf\circ\gamma^{-1}-g\circ\pi \]
is contained in $B_{\alpha\times c}$. 
Therefore we have $n[f]=n[f\circ\gamma^{-1}]$ in 
$K^0(X\times\Z_m,\alpha\times c)$. 
Since the $K^0$-group is torsion free, 
it follows that $[f]$ belongs to the kernel of $\id-\gamma^*$. 
By virtue of \cite[Lemma 3.6]{M2}, we get 
\[ T(K^0(X\times\Z_m,\alpha\times c)/\pi^*(K^0\xa))
=\Ker(\id-\gamma^*)/\pi^*(K^0\xa)\cong\Z_m. \]

Put $U=X\times \{0\}$. Then it is easy to see 
\[ f_0-f_0\circ\gamma^{-1}=
1_U-1_U\circ(\alpha\times c)^{-1}\in B_{\alpha\times c}. \]
Thus $[f_0]$ is in the kernel of $\id-\gamma^*$. 
By \cite[Lemma 3.6]{M2} and its proof, 
we can conclude that $[f_0]$ is the generator. 
\end{proof}

By using the lemma above, 
the computation of $PS(\alpha\times c)$ is carried out. 

\begin{lem}\label{PSofext}
Let $\xa$ be a Cantor minimal system and 
let $c:X\rightarrow\Z_2$ be a continuous map. 
Suppose $PS(\alpha)\neq PS(\alpha\times c)$. 
Then, there exists $n\in\N$ such that 
$2^{n-1}\in PS(\alpha)$, $2^n\notin PS(\alpha)$ 
and $PS(\alpha\times c)=2P(\alpha)\cup P(\alpha)$. 
Furthermore, when $2^{n-1}[f]$ equals $[1_X]$ in $K^0\xa$, 
we have $[c]=[f]+2K^0\xa$ in $K^0\xa/2K^0\xa$. 
\end{lem}
\begin{proof}
We follow the notation used in the lemma above. 
If $[c]$ is zero in $K^0\xa/2K^0\xa$, then 
obviously $PS(\alpha\times c)$ agrees with $PS(\alpha)$. 
Hence we may assume that $[c]$ is not zero, and that 
$\alpha\times c$ is minimal on $X\times\Z_2$. 

Suppose $p\in PS(\alpha\times c)\setminus PS(\alpha)$. 
Since $[1_{X\times\Z_2}]$ is divisible by $p$ 
in $K^0(X\times\Z_2,\alpha\times c)$, there exists 
$f\in C(X\times\Z_2,\Z)$ such that $p[f]=[1_{X\times\Z_2}]$. 
By $p\notin PS(\alpha)$, we have $[f]\notin\pi^*(K^0\xa)$. 
But $p[f]=[1_{X\times\Z_2}]=\pi^*([1_X])$, 
which means that $[f]$ gives a torsion element 
of $K^0(X\times\Z_2,\alpha\times c)/\pi^*(K^0\xa)$. 
Therefore, by the lemma above, we can see that 
$p$ is even and $2[f]\in\pi^*(K^0\xa)$. 
It follows that $p/2$ belongs to $PS(\alpha)$. 
Hence we can conclude that $PS(\alpha\times c)$ agrees 
with $2PS(\alpha)\cup PS(\alpha)$. 
Note that it also follows that $[f]$ is the generator of 
\[ T(K^0(X\times\Z_2,\alpha\times c)/
\pi^*(K^0\xa))\cong\Z_2. \]
If $[1_X]$ is infinitely 2-divisible in $K^0\xa$, then 
$PS(\alpha)$ contains $2PS(\alpha)$, which contradicts 
$PS(\alpha\times c)\neq PS(\alpha)$. 
Thus, there exists $n\in\N$ such that 
$2^{n-1}\in PS(\alpha)$ and $2^n\notin PS(\alpha)$. 

Suppose that $2^{n-1}[f]$ equals $[1_X]$ in $K^0\xa$. 
Since $2^n$ belongs to $PS(\alpha\times c)$, 
there exists $f_1\in C(X\times\Z_2,\Z)$ such that 
$2[f_1]=[f\pi]$ in $K^0(X\times\Z_2,\alpha\times c)$. 
As shown in the preceding paragraph, 
$[f_1]$ is the generator of 
\[ T(K^0(X\times\Z_2,\alpha\times c)/
\pi^*(K^0\xa))\cong\Z_2. \]
By the lemma above, 
\[ f_0(x,k)=\begin{cases}
1 & c(\alpha^{-1}(x))=1\mbox{ and }k=0 \\
0 & \mbox{otherwise} \end{cases} \]
is also the generator, and so we have $[f_1]-[f_0]\in\pi^*(K^0\xa)$. 
Observe that 
\[ 2[f_0]=[f_0]+[f_0\circ\gamma^{-1}]=[g_0\pi], \]
where $g_0\in C(X,\Z)$ is given by 
\[ g_0(x)=\begin{cases}
1 & c(x)=1 \\
0 & c(x)=0. \end{cases} \]
It follows that $2^{n-1}[g_0\pi]=2^n[f_0]$ is equal to 
$2^n[f_1]=2^{n-1}[f\pi]$ modulo $2^n\pi^*(K^0\xa)$. 
Hence $[g_0]-[f]$ belongs to $2K^0\xa$, because $\pi^*$ is injective. 
As a consequence, we get 
\[ [c]=[g_0]+2K^0\xa=[f]+2K^0\xa. \]
\end{proof}

\begin{lem}\label{combina}
Let $m_1,m_2,\dots,m_k$ be natural numbers and 
let $q$ be their greatest common divisor. 
Let $\chi_1,\chi_2,\dots,\chi_k$ be elements of $\Z_2$. 
Suppose that either of the following holds: 
\begin{enumerate}
\item There exist $i,j\in \{1,2,\dots,k\}$ such that 
$m_i\in2\N$, $\chi_i\neq0$ and $m_j\notin2\N$. 
\item There exist $i,j\in \{1,2,\dots,k\}$ such that 
$m_i\notin2\N$, $\chi_i=0$ and $\chi_j\neq0$. 
\end{enumerate}
Then, there exists a natural number $N\in\N$ such that, 
for every $n\geq N$ and $\chi\in\Z_2$, 
we can find $l_1,l_2,\dots,l_k\in\N$ so that 
\[ l_1m_1+l_2m_2+\dots+l_km_k=nq \]
and 
\[ l_1\chi_1+l_2\chi_2+\dots+l_k\chi_k=\chi. \]
\end{lem}
\begin{proof}
Let us consider the case (1). 
We may assume $m_1\in2\N$, $\chi_1\neq0$ and $m_2\notin2\N$. 
It is clear that there exists $N\in\N$ such that 
\[ \{nq:n\geq N\}\subset
\left\{\sum_{i=1}^kl_im_i:l_i\in\N, l_1>m_2\right\}. \]
Suppose that $n\geq N$ and $\chi$ are given. 
We can find $l_1,l_2,\dots,l_k\in\N$ so that 
\[ l_1m_1+l_2m_2+\dots+l_km_k=nq. \]
If $l_1\chi_1+l_2\chi_2+\dots+l_k\chi_k$ is equal to $\chi$, 
then we have nothing to do. 
Otherwise, by replacing $l_1$ with $l_1-m_2$ and 
$l_2$ with $l_2+m_1$, the assertion follows. 

In the case (2), we get the conclusion in a similar fashion. 
\end{proof}

Let 
\[ {\cal P}=\{X(v,k):v\in V,k=1,2,\dots,h(v)\} \]
be a Kakutani-Rohlin partition for a Cantor minimal system 
$\xa$. 
Let $\phi:X\rightarrow\Isom(\T)$ be a continuous map. 
Suppose that $o(\phi)$ is constant on each clopen set of ${\cal P}$. 
For every $v\in V$, 
\[ o(\phi)(x)+o(\phi)(\alpha(x))+\dots+o(\phi)(\alpha^{h(v)-1}(x)) \]
does not depend on $x\in X(v,1)$. 
We write this value by $o(\phi)_v\in\Z_2$. 

\begin{lem}\label{keyXT}
Let $\xa$ and $\yb$ be Cantor minimal systems and 
let $\phi:X\rightarrow\Isom(\T)$ and $\psi:Y\rightarrow\Isom(\T)$ 
be continuous maps. 
Suppose that 
\[ {\cal P}=\{X(v,k):v\in V,k=1,2,\dots,h(v)\} \]
and 
\[ {\cal Q}=\{Y(w,l):w\in W,l=1,2,\dots,h(w)\} \]
are Kakutani-Rohlin partitions for $\xa$ and $\yb$, and 
that $\pi:{\cal Q}\rightarrow{\cal P}$ is a bijection 
satisfying the following. 
\begin{enumerate}
\item $\pi(Y(w,l+1))=\alpha(\pi(Y(w,l)))$ for all $w\in W$ and 
$l=1,2,\dots,h(w)-1$. 
\item For all $v\in V$, we have $o(\phi)_v=\sum o(\phi)_w$, 
where the summation runs over all $w\in W$ such that 
$\pi(Y(w,1))$ belongs to the tower $v$. 
\end{enumerate}
Then there exist a homeomorphism $\sigma:X\rightarrow Y$ and 
a continuous map $\omega:X\rightarrow\Isom(\T)$ such that 
$\sigma(\pi(U))=U$ for all $U\in{\cal Q}$ and 
\[ \max_{x\in X}\max_{t\in\T}|\psi_{\sigma(x)}\omega_x(t)
-\omega_{\alpha(x)}\phi_x(t)|<\max_{v\in V}h(v)^{-1}. \]
\end{lem}
\begin{proof}
Choose a homeomorphism $\sigma:X\rightarrow Y$ so that 
$\sigma(\pi(U))=U$ for all $U\in{\cal Q}$. 

Let $v\in V$ and $x\in X(v,1)$. By (1) and (2), 
there exists $\kappa(x)\in\T$ such that 
\[ R_{\kappa(x)}=
\left(\psi_{\sigma(\alpha^{h(v)-1}(x))}
\psi_{\sigma(\alpha^{h(v)-2}(x))}\dots\psi_{\sigma(x)}\right)
\left(\phi_{\alpha^{h(v)-1}(x)}\phi_{\alpha^{h(v)-2}(x)}
\dots\phi_x\right)^{-1}. \]
Then $\kappa$ is a continuous function from $\alpha(R({\cal P}))$ to 
$\T$. 
Let $\tilde{\kappa}:\alpha(R({\cal P}))\rightarrow\R$ be 
a continuous function such that 
\[ \tilde{\kappa}(x)+\Z=\kappa(x) \ \mbox{ and } \ 
-1<\tilde{\kappa}(x)<1 \]
for all $x\in\alpha(R({\cal P}))$. 
Define $\chi(v,j)\in\Z_2$ by 
\[ \chi(v,j)=\sum_{i=j}^{h(v)-1}o(\psi_{\sigma(\alpha^i(x))}), \]
where $x$ belongs to $X(v,1)$ and $\chi(v,j)$ does not depend 
on the choice of $x$. 
Define $\eta:X\rightarrow\T$ by  
\[ \eta(\alpha^j(x))=-(-1)^{\chi(v,j)}
\frac{j}{h(v)}\tilde{\kappa}(x)+\Z \]
for all $x\in X(v,1)$ and $j=0,1,\dots,h(v)-1$. 
Then $\eta$ is a $\T$-valued continuous function on $X$. 
Since 
\[ o(\psi_{\sigma(\alpha^j(x))})+\chi(v,j)=\chi(v,j+1), \]
we have 
\[ \left|\psi_{\sigma(\alpha^j(x))}(R_{\eta(\alpha^j(x))}(t))-
R_{\eta(\alpha^{j+1}(x))}(\psi_{\sigma(\alpha^j(x))}(t))\right|
\leq\frac{1}{h(v)}|\tilde{\kappa}(x)|<\frac{1}{h(v)} \]
for all $t\in\T$, $x\in X(v,1)$ and $j=0,1,\dots,h(v)-1$. 
Define a continuous map $\omega:X\rightarrow\Isom(\T)$ by 
\[ \omega_x=\id \]
for all $x\in\alpha(R({\cal P}))$ and 
\[ \omega_{\alpha^j(x)}=R_{\eta(\alpha^j(x))}
\left(\psi_{\sigma(\alpha^{j-1}(x))}\psi_{\sigma(\alpha^{j-2}(x))}
\dots\psi_{\sigma(x)}\right)
\left(\phi_{\alpha^{j-1}(x)}\phi_{\alpha^{j-2}(x)}
\dots\phi_x\right)^{-1} \]
for all $x\in X(v,1)$ and $j=1,2,\dots,h(v)-1$. 
Then, one checks that 
\[ \left|\psi_{\sigma(\alpha^j(x))}\omega_{\alpha^j(x)}(t)
-\omega_{\alpha^{j+1}(x)}\phi_{\alpha^j(x)}(t)\right|
<\frac{1}{h(v)} \]
for all $t\in\T$, $x\in X(v,1)$ and $j=0,1,\dots,h(v)-1$. 
\end{proof}

The following is the main theorem of this section. 

\begin{thm}\label{wacXT}
Let $\xa$ and $\yb$ be Cantor minimal systems and 
let $\phi:X\rightarrow\Isom(\T)$ and $\psi:Y\rightarrow\Isom(\T)$ 
be continuous maps. 
The following are equivalent. 
\begin{enumerate}
\item There exist homeomorphisms $\sigma_n$ from $X\times\T$ to 
$Y\times\T$ such that $\sigma_n(\alpha\times\phi)\sigma_n^{-1}$ 
converges to $\beta\times\psi$ in $\Homeo(Y\times\T)$.
\item $PS(\beta)$ is contained in $PS(\alpha)$, and 
either of the following conditions is satisfied: 
\begin{enumerate}
\item $[o(\psi)]_\beta=0$ and $[o(\phi)]_\alpha=0$
\item $[o(\psi)]_\beta\neq0$ and there exist $n\in\N$, 
$f\in C(X,\Z)$ and $g\in C(Y,\Z)$ such that 
\[ 2^{n-1}[f]_\alpha=[1_X]_\alpha, \ 
[f]_\alpha+2K^0\xa=[o(\phi)]_\alpha \]
and 
\[ 2^{n-1}[g]_\beta=[1_Y]_\beta, \ 
[g]_\beta+2K^0\yb=[o(\psi)]_\beta. \]
\item $[o(\psi)]_\beta\neq0$ and 
$PS(\beta\times o(\psi))=PS(\beta)$. 
\end{enumerate}
\end{enumerate}
\end{thm}
\begin{proof}
(1)$\Rightarrow$(2). 
Note that every homeomorphism from $X\times\T$ to $Y\times\T$ is 
of the form $\sigma\times\omega$, 
where $\sigma:X\rightarrow Y$ is a homeomorphism and 
$\omega$ is a continuous map from $X$ to $\Homeo(\T)$. 

By the same proof as (1)$\Rightarrow$(2) of Theorem \ref{wac}, 
we see that $PS(\beta)=PS(\beta\times\psi)$ is contained 
in $PS(\alpha)=PS(\alpha\times\phi)$. 

At first, suppose that $[o(\psi)]_\beta$ is not zero 
and there exist $n\in\N$ and $g\in C(Y,\Z)$ such that 
\[ 2^{n-1}[g]_\beta=[1_Y]_\beta, \ 
[g]_\beta+2K^0\yb=[o(\psi)]_\beta. \]
We will show (b). 
There exists a Kakutani-Rohlin partition 
\[ {\cal Q}=\{Y(w,l):w\in W,l=1,2,\dots,h(w)\} \]
such that the following hold: 
\begin{itemize}
\item For all $w\in W$, $h(w)$ is divisible by $2^{n-1}$. 
Put $m_w=2^{1-n}h(w)$. 
\item The function $o(\psi)$ is constant 
on each clopen set belonging to ${\cal Q}$. 
\item $o(\psi)_w$ is equal to $m_w+2\Z$ for all $w\in W$. 
\end{itemize}
By the assumption, 
we can find a homeomorphism $\sigma:X\rightarrow Y$ and 
a continuous map $\omega:X\rightarrow\Homeo(\T)$ such that 
\[ \sigma\alpha\sigma^{-1}(U)=\beta(U) \]
for all $U\in {\cal Q}$ and 
\[ |\omega_{\alpha(x)}(\phi_x(t))
-\psi_{\sigma(x)}(\omega_x(t))|<\frac{1}{2} \]
for all $(x,t)\in X\times\T$. 
Hence we have 
\[ o(\omega)(\alpha(x))+o(\phi)(x)
=o(\psi)(\sigma(x))+o(\omega)(x) \]
for all $x\in X$, 
which implies that $[o(\phi)]_\alpha=[o(\psi\sigma)]_\alpha$. 
Defining $X(w,l)=\sigma^{-1}(Y(w,l))$, 
we obtain a Kakutani-Rohlin partition 
\[ {\cal P}=\{X(w,l):x\in W,l=1,2,\dots,h(w)\} \]
for $\xa$. 
Put 
\[ f(x)=\begin{cases}
1 & x\in X(w,l),l\in\{1,2,\dots,m_w\} \\
0 & \mbox{otherwise}. \end{cases} \]
Then $2^{n-1}[f]_\alpha=[1_X]_\alpha$, and 
\[ [o(\phi)]_\alpha=[o(\psi\sigma)]_\alpha
=\sum_{w\in W}m_w[1_{X(w,1)}]_\alpha+2K^0\xa=[f]_\alpha+2K^0\xa, \]
which means that (b) holds. 

If $[o(\psi)]_\beta$ is zero, then we can find 
a Kakutani-Rohlin partition ${\cal Q}$ such that 
$o(\psi)_w=0$ for all $w\in W$. 
The same argument as above implies $[o(\phi)]_\alpha=0$. 
Thus (a) holds. 

In the other case, (c) is immediate from Lemma \ref{PSofext}. 

(2)$\Rightarrow$(1). 
Suppose that $PS(\beta)$ is contained in $PS(\alpha)$. 
We will modify the argument of Theorem \ref{wac} 
and apply Lemma \ref{keyXT}. 
Take a clopen partition ${\cal F}$ of $Y$ arbitrarily and 
let $\varepsilon>0$. 
It suffices to find a homeomorphism $\sigma:X\rightarrow Y$ and 
a continuous map $\omega:X\rightarrow\Isom(\T)$ such that 
\[ \sigma\alpha\sigma^{-1}(U)=\beta(U) \]
for all $U\in{\cal F}$ and 
\[ |\psi_{\sigma(x)}(\omega_x(t))
-\omega_{\alpha(x)}(\phi_x(t))|<\varepsilon \]
for all $(x,t)\in X\times\T$. 

Let us assume (b). 
There exists a Kakutani-Rohlin partition 
\[ {\cal Q}=\{Y(w,k):w\in W,l=1,2,\dots,h(w)\} \]
such that the following conditions are satisfied. 
\begin{itemize}
\item For all $w\in W$, $h(w)$ is divisible by $2^{n-1}$. 
Put $m_w=2^{1-n}h(w)$. 
\item The function $o(\psi)$ is constant 
on each clopen set belonging to ${\cal Q}$. 
\item $o(\psi)_w$ is equal to $m_w+2\Z$ for all $w\in W$. 
\item $\widetilde{\cal Q}$ is finer than ${\cal F}$, where 
\[ \widetilde{\cal Q}=\{Y(w,k):w\in W,l=1,2,\dots,h(w)-1\}
\cup\{R({\cal Q})\}. \]
\end{itemize}
Let $p$ be the greatest common divisor of $m_w$'s. 
Let 
\[ {\cal P}=\{X(v,k):v\in V,k=1,2,\dots,h(v)\} \]
be a Kakutani-Rohlin partition for $\xa$ 
which satisfies the following. 
\begin{itemize}
\item For all $v\in V$, $h(v)$ is divisible by $2^{n-1}p$ 
and sufficiently large. 
\item The function $o(\phi)$ is constant 
on each clopen set belonging to ${\cal P}$. 
\item $o(\phi)_v$ is equal to $2^{1-n}h(v)+2\Z$ for all $v\in V$. 
\end{itemize}
In the same way as in Theorem \ref{wac}, 
we can find natural numbers $a_{v,w}$ so that 
\[ h(v)=\sum_{w\in W}a_{v,w}h(w). \]
It follows that 
\[ o(\phi)_v=2^{1-n}h(v)+2\Z
=2^{1-n}\sum_{w\in W}a_{v,w}h(w)+2\Z
=\sum_{w\in W}a_{v,w}o(\psi)_w. \]
Therefore Lemma \ref{keyXT} applies and 
the required $\sigma:X\rightarrow Y$ and 
$\omega:X\rightarrow\Isom(\T)$ are obtained. 

When we assume (a), a similar argument is valid 
(or we can say that it follows from (2)$\Rightarrow$(1) of 
Theorem \ref{proto}). 

Finally, let us consider the case (c). 
Let ${\cal Q}$ be a Kakutani-Rohlin partition for $\yb$ such that 
$o(\psi)$ is constant on each clopen set belonging to ${\cal Q}$ 
and $\widetilde{\cal Q}$ is finer than ${\cal F}$. 
Let $p$ be the greatest common divisor of $h(w)$'s and 
let $p=2^{n-1}q$ with $n\in\N$ and $q\notin2\N$. 
Put $m_w=2^{1-n}h(w)$. 
Note that the greatest common divisor of $m_w$'s is $q$. 
If $o(\psi)_w=m_w+2\Z$ for all $w\in W$, then 
there exists $f\in C(Y,\Z)$ such that 
\[ 2^{n-1}[f]_\beta=[1_Y]_\beta \]
and 
\[ [f]_\beta+2K^0\yb=[o(\psi)]_\beta. \]
By Lemma \ref{torsion} and \ref{PSofext}, we have 
$PS(\beta\times o(\psi))\neq PS(\beta)$, which is a contradiction. 
Hence there exists $w\in W$ such that $o(\psi)_w\neq m_w+2\Z$. 
Moreover, there exists $w\in W$ such that 
$o(\psi)_w$ is not zero, because $[o(\psi)]_\beta$ is not zero. 
It follows that we have either of the following. 
\begin{itemize}
\item There are $w_1,w_2\in W$ such that $m_{w_1}$ is even, 
$o(\psi)_{w_1}\neq 0$ and $m_{w_2}$ is odd. 
\item There are $w_1,w_2\in W$ such that $m_{w_1}$ is odd, 
$o(\psi)_{w_1}=0$ and $o(\psi)_{w_2}\neq 0$. 
\end{itemize}
Then Lemma \ref{combina} applies and yields a natural number $N$. 
Let 
\[ {\cal P}=\{X(v,k):v\in V,k=1,2,\dots,h(v)\} \]
be a Kakutani-Rohlin partition for $\xa$ 
which satisfies the following. 
\begin{itemize}
\item For all $v\in V$, $h(v)$ is divisible by $p$ and 
greater than $Np$. 
\item The function $o(\phi)$ is constant on each clopen set 
belonging to ${\cal P}$. 
\end{itemize}
By the choice of $N$, there exist natural numbers $a_{v,w}$ such that 
\[ h(v)=\sum_{w\in W}a_{v,w}h(w) \]
and 
\[ o(\phi)_v=\sum_{w\in W}a_{v,w}o(\psi)_w \]
for every $v\in V$. 
Therefore we can apply Lemma \ref{keyXT} and 
complete the proof. 
\end{proof}

Combining Lemma \ref{PSofext} and Theorem \ref{wacXT}, 
we get the following. 

\begin{cor}
Let $\xa$ and $\yb$ be Cantor minimal systems and 
let $\phi:X\rightarrow\Isom(\T)$ and $\psi:Y\rightarrow\Isom(\T)$ 
be continuous maps. 
Then $\alpha\times\phi$ and $\beta\times\psi$ are 
weakly approximately conjugate if and only if 
either of the following holds. 
\begin{enumerate}
\item $PS(\alpha)=PS(\beta)$ and 
both $\alpha\times\phi$ and $\beta\times\psi$ are orientation preserving. 
\item $PS(\alpha)=PS(\beta)$, 
$PS(\alpha\times o(\phi))=PS(\beta\times o(\psi))$ and 
neither $\alpha\times\phi$ nor $\beta\times\psi$ are 
orientation preserving. 
\end{enumerate}
\end{cor}

In the proof of Lemma \ref{keyXT}, it was essentially used that 
the cocycles $\phi$ and $\psi$ are $\Isom(\T)$-valued. 
Unfortunately, our proof does not work 
for cocycles with values in $\Homeo(\T)$. 
We do not know when two minimal homeomorphisms on $X\times\T$ are 
weakly approximately conjugate in general.

\flushleft{
\textit{e-mail: matui@math.s.chiba-u.ac.jp \\
Graduate School of Science and Technology,\\
Chiba University,\\
1-33 Yayoi-cho, Inage-ku,\\
Chiba 263-8522,\\
Japan. }}

\end{document}